\documentclass[reqno]{article}
\usepackage{amsmath,amsfonts,amssymb,amsthm}

\newtheorem{theorem}{Theorem}

\theoremstyle{definition}
\newtheorem{remark}{Remark}

\newcommand\eps{\varepsilon}
\renewcommand\le{\leqslant}
\renewcommand\ge{\geqslant}
\newcommand\cc{{\mathrm{c}}}
\newcommand\E{{\mathbb E}}
\renewcommand\Pr{{\mathbb P}}

\newcommand\cF{\mathcal{F}}

\newcommand\cI{{\mathcal I}}
\newcommand\cJ{{\mathcal J}}

\newcommand\cS{{\mathcal S}}

\begin{document}
\title{The Janson inequalities for general up-sets}
\author{Oliver Riordan%
\thanks{Mathematical Institute, University of Oxford, 24--29 St Giles', Oxford OX1 3LB, UK.
E-mail: {\tt riordan@maths.ox.ac.uk}.}
\ and Lutz Warnke%
\thanks{Department of Pure Mathematics and Mathematical Statistics,
Wilberforce Road, Cambridge CB3 0WB, UK.
E-mail: {\tt L.Warnke@dpmms.cam.ac.uk}.}}
\date{March 5, 2012; revised April 17, 2013}
\maketitle

\begin{abstract}
Janson and Janson, {\L}uczak and Ruci\'nski proved several inequalities 
for the lower tail of the distribution of the number of events that hold, when all the
events are up-sets (increasing events) of a special form --
each event is the intersection of some subset of a single set of independent
events (i.e., a principal up-set). We show that these 
inequalities in fact hold for arbitrary up-sets, by modifying existing proofs
to use only positive correlation, avoiding the need to assume positive correlation conditioned
on one of the events.
\end{abstract}

Let $(\Omega,\cF,\Pr)$ be a probability space, and $\cI\subseteq\cF$ a collection of events
with the following properties:
\begin{equation}\label{Icorr}
 A,B\in \cI \implies \Pr(A\cap B)\ge \Pr(A)\Pr(B)
\end{equation}
and
\begin{equation}\label{Iclosed}
 A,B\in \cI \implies  A\cap B\in \cI \hbox{ and } A\cup B\in \cI.
\end{equation}
The standard and most important example is when $\Omega=\{0,1\}^S$
is a product probability space (with product measure), and $\cI$ is the collection
of all increasing events, i.e., events $A$ such that $\omega\in A$ and $\omega\le\omega'$ pointwise
imply $\omega'\in A$. (Of course, one can instead take all decreasing events.) 
Then \eqref{Icorr} is simply Harris's Lemma~\cite{Harris} (also known as
Kleitman's Lemma~\cite{Kleitman}). There are other examples, such as increasing events in random cluster
models with parameter $q\ge 1$; see~\cite{FKG}.

Let $A_1,\ldots,A_k\in \cI$, write $I_i$ for the indicator function of $A_i$, and set
\[
 X=\sum_i I_i, \quad\quad \mu=\E(X)=\sum_i \Pr(A_i)
\]
and
\[
 \Delta=\sum_i\sum_{j\sim i} \Pr(A_i\cap A_j),
\]
where we write $i\sim j$ if $i\ne j$ and $A_i$ and $A_j$ are dependent. (Note that we sum over \emph{ordered} pairs,
and exclude the term $i=j$. These conventions are not universal!)

\begin{theorem}\label{th1}
Under the conditions above we have
\begin{equation}\label{i1}
 \Pr(X=0) \le \exp(-\mu+\Delta/2),
\end{equation}
and, for any $0\le t\le \mu$,
\begin{equation}\label{i2}
 \Pr(X\le \mu-t) \le \exp\left(-\frac{\phi(-t/\mu)\mu^2}{\mu+\Delta}\right) \le \exp\left(-\frac{t^2}{2(\mu+\Delta)}\right),
\end{equation}
where $\phi(x)=(1+x)\log(1+x)-x$ with $\phi(-1)=1$.
\end{theorem}

When the events $A_i$ are \emph{principal} up-sets, i.e.,
events in a product space $\{0,1\}^S$
of the form $A_i=\{\omega: \omega_x=1 \text{ for all }x\in \alpha_i\}$ for some $\alpha_i\subseteq S$,
the first inequality in \eqref{i2} is the well known Janson inequality~\cite{Janson}.
The second is also given in~\cite{Janson}; for other
convenient weaker forms see~\cite{JLR_book}. Under the same assumptions, \eqref{i1} was proved earlier by
Janson, {\L}uczak and Ruci\'nski~\cite{JLR_ineq}.
We shall prove Theorem~\ref{th1} by modifying the proofs of these inequalities to avoid applying Harris's Lemma
to the conditional measure $\Pr(\cdot\mid A_i)$.

\begin{proof}
We begin with a simple observation that, in the standard setting, follows from the equality conditions in Harris's Lemma.
Indeed, we claim that, for each~$i$,
\begin{equation}\label{claim}
  A_i\hbox{ is independent of the set }\cS_i=\{A_j: j\ne i,\, j\not\sim i\} .
\end{equation}
To see this, note first that if $A,B,C\in \cI$ and $A$ is independent of $B$ and of $C$, then
\begin{equation}\label{e1}
 \Pr(A\cap(B\cap C))+\Pr(A\cap(B\cup C)) = \Pr(A\cap B)+\Pr(A\cap C) = \Pr(A)(\Pr(B)+\Pr(C)).
\end{equation}
Since, by \eqref{Iclosed}, $B\cap C$ and $B\cup C$ are in $\cI$, by \eqref{Icorr} we have
\[
 \Pr(A\cap(B\cap C))\ge\Pr(A)\Pr(B\cap C) \hbox{ and } \Pr(A\cap (B\cup C))\ge \Pr(A)\Pr(B\cup C).
\]
The sum of these inequalities is the equality \eqref{e1}, so both inequalities are equalities,
and in particular $A$ is independent of $B\cap C$. Since $A_i$ is independent of each $A_j\in \cS_i$
by definition, it follows inductively that $A_i$ is independent of any intersection of events $A_j\in \cS_i$,
so $A_i$ is independent of the set $\cS_i$ of events, as claimed.

To prove \eqref{i1} we modify the argument given by Boppana and Spencer~\cite{BS}, following the presentation
in~\cite{AS}. Following Boppana and Spencer, set
$r_i=\Pr(A_i\mid A_1^\cc\cap\cdots\cap A_{i-1}^\cc)$, so $\Pr(X=0)=\prod_{i=1}^k (1-r_i)\le \exp(-\sum r_i)$.
It suffices to show that for each $i$ we have
\begin{equation}\label{aim}
 r_i\ge \Pr(A_i) -\sum_{j<i,\,j\sim i} \Pr(A_i\cap A_j),
\end{equation}
since the sum of the final expression is exactly $\mu-\Delta/2$.
Fix $i$. Still following~\cite{BS,AS}, set
\[
D=\bigcap_{j<i,\,j\sim i} A_j^\cc \hbox{\quad and\quad} R = \bigcap_{j<i,\, j\not\sim i} A_j^\cc,
\]
noting that, by \eqref{claim}, $A_i$ is independent of $R$. Then
\begin{equation}\label{e2}
 r_i=\Pr(A_i\mid D\cap R) = \frac{\Pr(A_i\cap D\cap R)}{\Pr(D\cap R)} 
 \ge \frac{\Pr(A_i\cap D\cap R)}{\Pr(R)} = \Pr(A_i\cap D\mid R).
\end{equation}
(We may assume that $\Pr(D\cap R)>0$, since otherwise $\Pr(X=0)=0$.)
At this point the original argument involves writing the final probability as $\Pr(A_i\mid R)\Pr(D\mid A_i\cap R)$.
Instead, we simply write
\begin{equation}\label{e3}
  \Pr(A_i\cap D\mid R) = \Pr(A_i\mid R) - \Pr(A_i\cap D^\cc \mid R).
\end{equation}
By \eqref{Iclosed}, $D^\cc=\bigcup_{j<i,\, j\sim i} A_j \in \cI$, 
so $A_i\cap D^\cc\in \cI$. Since $R^\cc\in \cI$,  
using \eqref{Icorr} and the union bound it follows that 
\begin{equation}\label{e4}
 \Pr(A_i\cap D^\cc \mid R) \le \Pr(A_i\cap D^\cc) = \Pr\left(A_i\cap \bigcup_{j<i,\,j\sim i} A_j\right)
 \le \sum_{j<i,\,j\sim i} \Pr(A_i\cap A_j).
\end{equation}
Recalling that $A_i$ is independent of $R$, combining \eqref{e2}, \eqref{e3} and \eqref{e4} 
gives \eqref{aim}, completing the proof of \eqref{i1}.

Turning to \eqref{i2}, fix $1\le i\le k$ and let
\[
 Y_i = I_i+\sum_{j\sim i} I_j \hbox{\quad and\quad} Z_i=\sum_{j\ne i,\,j\not\sim i} I_j,
\]
so $X=Y_i+Z_i$, with $Z_i$ containing the terms independent of $I_i$ and $Y_i$ the others (including $I_i$ itself).
In the proof of \eqref{i2} given in~\cite{JLR_book},
the only step in which anything is assumed about the events $A_i$ is (2.20) on page 32,
where it is shown (in our notation) that for $s\ge 0$ and each $1\le i\le k$,
\begin{equation}\label{aim2}
 \E(I_i e^{-sX}) \ge \E(I_i e^{-sY_i}) \E(e^{-sX}).
\end{equation}
Proceeding much as in the proof of \eqref{aim}, note that
\begin{equation}\label{e5}
 \frac{ \E(I_i e^{-sX}) }{\E(e^{-sX})} =
  \frac{ \E(I_i e^{-sX}) }{\E(e^{-sY_i}e^{-sZ_i})}  \ge
  \frac{ \E(I_i e^{-sX}) }{\E(e^{-sZ_i})}.
\end{equation}
Also,
\[
 I_ie^{-sX} = I_ie^{-sY_i}e^{-sZ_i} = I_ie^{-sZ_i} - I_ie^{-sZ_i}(1 - e^{-sY_i}) = I_ig -fg,
\]
where
\[ 
 f=I_i(1 - e^{-sY_i}) \hbox{\quad and \quad} g=e^{-sZ_i}.
\]
Now from \eqref{claim}, $I_i$ and $Z_i$ are independent, so $\E(I_ig)=\E(I_i)\E(g)$.
We may write $f$ in the form $f=v_0+\sum_j (v_j-v_{j-1})J_j$, where 
$0 \le v_0<v_1<\cdots$ are the distinct values taken by $f$, and each $J_j$ 
is the indicator function of the event $B_j=\{f\ge v_j\}$. 
Note that any such $B_j$ may be expressed as 
$\bigcup_{\alpha \in \cJ} \bigcap_{i \in \alpha} A_i$ for some set $\cJ$
of subsets of $\{1,2,\ldots,k\}$, so 
\eqref{Iclosed} implies $B_j \in \cI$. 
Writing $1-g$ in an analogous form, it follows from \eqref{Icorr} that 
$\E(f(1-g))\ge \E(f)\E(1-g)$, so $\E(fg) \le  \E(f)\E(g)$. 
Hence, 
\[
 \E(I_ie^{-sX}) = \E(I_ig-fg) \ge \E(I_i)\E(g)-\E(f)\E(g).
\]
Using \eqref{e5} for the first step this gives
\[
  \frac{ \E(I_i e^{-sX}) }{\E(e^{-sX})} \ge  \frac{ \E(I_i e^{-sX}) }{\E(g)} \ge \E(I_i)-\E(f)
 = \E(I_i-f) = \E(I_ie^{-sY_i}).
\]
This is exactly \eqref{aim2}, and the rest of the proof in \cite{JLR_book} is unaltered.
\end{proof}

\begin{remark}
We have stated two of the best-known and cleanest forms of the inequalities
in Theorem~\ref{th1}. Let us note that other forms also hold in the present more general
context. For example, the inequality \eqref{aim} leads to the bound
\begin{equation}\label{i1a}
 \Pr(X=0) \le \prod_{i=1}^k(1-\Pr(A_i))\exp\left(\tfrac{1}{1-\eps}\tfrac{\Delta}{2}\right),
\end{equation}
where $\eps=\max_i \Pr(A_i)$. This bound was given by Boppana and Spencer~\cite{BS} (for $\eps=1/2$); see
also~\cite{AS,JLR_book}.

Furthermore, \eqref{aim2} is the only step in the proof of Lemma 1 of~\cite{JLR_ineq} that requires
any assumptions about the $A_i$. 
Hence this result, which is slightly stronger than \eqref{i1}, also holds in the present setting,
giving (in our notation)
\begin{equation}\label{i2a}
 \log\Pr(X=0) \le -\sum_i \E\left(\frac{I_i}{I_i+\sum_{j\sim i} I_j}\right).
\end{equation} 
\end{remark}

\begin{remark}
A key feature of the various Janson inequalities is that when $\Delta$ is small, then
the bounds are close to best possible, since $\Pr(X=0)\ge \prod_i (1-\Pr(A_i))$.
We have not stressed this since it is well known that this lower bound applies
to general up-sets $A_i$, by Harris's Lemma. Similarly, it applies whenever
the $A_i$ are in some collection $\cI$ of events satisfying \eqref{Icorr} and \eqref{Iclosed}.
\end{remark}

\begin{remark}
The bounds in \eqref{i2} can be extended to the weighted case $X=\sum_{i}c_iI_i$ with 
positive $c_i$, studied, e.g., in~\cite{KV}: we obtain
\begin{equation}\label{i2e}
 \Pr(X\le \mu-t) \le \exp\left(-\frac{\phi(-t/\mu)\mu^2}{\overline{\Delta}}\right)
 \le \exp\left(-\frac{t^2}{2\overline{\Delta}}\right),
\end{equation}
where $\mu=\E(X)$ and
\[
 \overline{\Delta}=\sum_i c_i^2\Pr(A_i) + \sum_i\sum_{j\sim i} c_ic_j\Pr(A_i\cap A_j)=
 \sum_i \E(J_i^2) + \sum_i\sum_{j\sim i} \E(J_iJ_j),
\]
for $J_\ell=c_\ell I_\ell$. (In applications, it may be convenient to note that $\sum_i \E(J_i^2) \le C \mu$ 
when $\max_i c_i \le C$.) 
The proof of \eqref{i2e} is a straightforward 
modification of that of Theorem~2.14 in~\cite{JLR_book}. The key inequality
is again \eqref{aim2}, now with $I_\ell$ replaced by $J_\ell$ in the definitions of
$X$, $Y_i$ and $Z_i$. The proof above carries over since all $c_i$ are positive.
Finally, in this setting \eqref{i2a} also holds, with $I_\ell$ replaced 
by $J_\ell$.  
\end{remark}

\end{document}